\documentclass[11pt,a4paper]{article}
\usepackage{latexsym}
\begin{document}
\newtheorem{thm}{Th{\'e}or{\`e}me}[section]
\newtheorem{cor}{Corollaire}[section]
\newtheorem{prop}{Proposition}[section]
\newtheorem{defi}{D{\'e}finition}[section]
\newtheorem{lem}{Lemme}[section]
\newtheorem{ex}{Exemple}[section]
\title{Quasi-big{\`e}bres de Lie et cohomologie d'alg\`{e}bre de Lie}
\author{Momo BANGOURA\footnote{Senior Associate of the Abdus Salam ICTP, Trieste, Italy}\\
D{\'e}partement de Math{\'e}matiques, Universit{\'e} de Conakry\\ BP 1147,
R{\'e}publique de Guin{\'e}e\footnote{Permanent address, E-mail: bangm59@yahoo.fr}}
\date{}
\maketitle
\begin{abstract}
Lie quasi-bialgebras are natural generalisations of Lie bialgebras introduced by Drinfeld. To any Lie quasi-bialgebra structure of finite-dimensional $(\mathcal{G}, \mu, \gamma, \phi)$,  correspond one Lie algebra structure on $\mathcal{D}= \mathcal{G}\oplus \mathcal{G^{*}}$, called the double of the given Lie quasi-bialgebra. We show that there exist on  $\Lambda\mathcal{G}$, the exterior algebra of $\mathcal{G}$, a $\mathcal{D}$-module structure and we establish an isomorphism of $\mathcal{D}$-modules between $\Lambda\mathcal{D}$ and $End(\Lambda\mathcal{G})$, $\mathcal{D}$ acting on $\Lambda\mathcal{D}$ by the adjoint action.
\begin{center}{\bf R{\'e}sum{\'e}}\end{center}
Les quasi-big{\`e}bres de Lie sont des g\'{e}n\'{e}ralisations naturelles, introduites par Drinfeld, des big\`{e}bres de Lie. A toute structure de quasi-big{\`e}bre de Lie $(\mathcal{G}, \mu, \gamma, \phi)$ de dimension finie, il correspond une structure
d'alg{\`e}bre de Lie  sur $\mathcal{D}= \mathcal{G}\oplus \mathcal{G^{*}}$, appel\'{e}e le double de la quasi-big{\`e}bre de Lie donn\'{e}e. On montre qu'il existe sur $\Lambda\mathcal{G}$, l'alg{\`e}bre ext{\'e}rieure de $\mathcal{G}$, une structure de $\mathcal{D}$-module et nous \'{e}tablissons un isomorphisme de  $\mathcal{D}$-modules entre $\Lambda\mathcal{D}$ et $End(\Lambda\mathcal{G})$, $\mathcal{D}$ agissant sur $\Lambda\mathcal{D}$  par l'action adjointe.
\end{abstract}
\newpage
\tableofcontents
\newpage
\section{Introduction}
Le but de ce travail est de construire pour une quasi-big{\`e}bre de Lie donn\'{e}e $(\mathcal{G}, \mu, \gamma, \phi)$, une repr\'{e}sentation naturelle de son double $\mathcal{D}$ sur l'alg{\`e}bre ext{\'e}rieure de $\mathcal{G}$, ou encore une structure de $\mathcal{D}$-module sur $\Lambda\mathcal{G}$, de telle sorte qu'il existe un isomorphisme de $\mathcal{D}$-modules entre $\Lambda\mathcal{D}$ et $End(\Lambda\mathcal{G})$, $\mathcal{D}$ agissant sur $\Lambda\mathcal{D}$  par l'action adjointe.

Les quasi-big{\`e}bres de Lie (\cite{Dr2}) ou quasi-big\`{e}bres jacobiennes (\cite{Ba1}, \cite{BK}, \cite{KS}) sont des g\'{e}n\'{e}ralisations naturelles des big\`{e}bres de Lie (\cite{Dr1}), introduites par Drinfeld comme \'{e}tant les limites classiques des alg\`{e}bres quasi-Hopf; contrairement aux big\`{e}bres de Lie, elles sont caract\'{e}ris\'{e}es par l'existence d'un d\'{e}faut d'identit\'{e} de co-Jacobi pour le co-crochet, qui est en fait le cobord d'un certain \'{e}l\'{e}ment de $\Lambda^{3}\mathcal{G}$, o\`{u} $\mathcal{G}$ est l'espace vectoriel sur lequel est d\'{e}finie la structure de quasi-big{\`e}bre de Lie, alors que pour les big\`{e}bres de Lie, ce d\'{e}faut est nul.

Dans la section 2, nous faisons un bref rappel de quelques notions fondamentales qui sont les outils de travail dans toute la suite, notamment le crochet de Schouten alg\'{e}brique, la cohomologie d'alg\`{e}bre de Lie.

Dans la section 3, nous rappelons la d\'{e}finition et les propri\'{e}t\'{e}s des quasi-big{\`e}bres de Lie et \`{a} partir d'une structure de quasi-big\`{e}bre de Lie donn\'{e}e $(\mathcal{G}, \mu, \gamma, \phi)$, nous d\'{e}finissons des op\'{e}rateurs de cohomologie sur $\Lambda\mathcal{G}$ et $\Lambda\mathcal{G}^{*}$, qui sont li\'{e}s par un ensemble de relations, cons\'{e}quences des axiomes de la structure de quasi-big\`{e}bre de Lie. Enfin, nous d\'{e}finissons le laplacien d'une quasi-big\`{e}bre de Lie, qui est une d\'{e}rivation de degr\'{e} 0 de $(\Lambda\mathcal{G}, \wedge)$ et de $(\Lambda\mathcal{G}, [ , ]^{\mu})$, o\`{u} $[ , ]^{\mu}$ est le crochet de Schouten alg\'{e}brique (\cite{KS}, \cite{KM}, \cite{Kz2})  d\'{e}fini \`{a} partir de la structure d'alg\`{e}bre de Lie sur $\mathcal{G}$; on montre qu'il commute avec l'op\'{e}rateur d'homologie de Chevalley-Eilenberg (\`{a} coefficients triviaux) d\'{e}fini \'{e}galement \`{a} partir de la structure d'alg\`{e}bre de Lie sur $\mathcal{G}$.

La section 4 recouvre l'essentiel du travail, \`{a} savoir la d\'{e}fintion d'une repr\'{e}sentation canonique de l'alg\`{e}bre de Lie double $\mathcal{D}$ d'une quasi-big\`{e}bre de Lie $(\mathcal{G}, \mu, \gamma, \phi)$ sur son alg\`{e}bre ext\'{e}rieure $\Lambda\mathcal{G}$ et l'\'{e}tablissement d'un isomorphisme de $\mathcal{D}$-modules entre $\Lambda\mathcal{D}$ et $End(\Lambda\mathcal{G})$, $\mathcal{D}$ agissant sur $\Lambda\mathcal{D}$  par l'action adjointe. Pour cela nous utilisons les constructions de (\cite{Lu}) bas\'{e}es sur la th\'{e}orie des alg\`{e}bres de Clifford (\cite{K-S}).

Dans toute la suite nous supposerons les structures d'alg\`{e}bre de Lie de dimension finie. Ainsi, si $(\mathcal{G}, \mu)$ est une alg\`{e}bre de Lie et  $\mathcal{G}^{*}$ son espace vectoriel dual, le crochet de dualit\'{e} entre $\Lambda\mathcal{G}$ et $\Lambda\mathcal{G}^{*}$ \'{e}tendant celui entre $\mathcal{G}$ et $\mathcal{G}^{*}$ est d\'{e}fini par
$$<\xi_{1}\wedge\xi_{2}\wedge...\wedge\xi_{m}, x_{1}\wedge x_{2}\wedge...\wedge x_{n}> = \delta_{m}^{n}det(<\xi_{i}, x_{j}>),$$
$\xi_{i}\in \mathcal{G}^{*}, i = 1,...,m, x_{j}\in\mathcal{G}, j = 1,...,n$.\\
Pour tous $X\in \Lambda\mathcal{G}$, notons par $\varepsilon_{X}\in End(\Lambda\mathcal{G})$ l'application d\'{e}finie par
$$Y\in \Lambda\mathcal{G} \rightarrow X\wedge Y \in \Lambda\mathcal{G},$$
et par $i_{X}\in End(\Lambda\mathcal{G}^{*})$ sa transpos\'{e}e d\'{e}finie par
$$<i_{X}A, Y> = <A, X\wedge Y>, \forall Y\in \mathcal{G}, \forall A\in \mathcal{G}^{*}.$$
\section{Pr\'{e}liminaires}
Dans cette section, nous rappelons certaines notions standard utiles pour la suite du travail.
\subsection{Crochet de Schouten alg\'{e}brique}
Soit $(\mathcal{G}, \mu)$ une alg\`{e}bre de Lie sur le corps {\it K}, suppos\'{e} \'{e}gal \`{a} {\it R} ou {\it C}, o\`{u} $\mu : \Lambda^{2}\mathcal{G}\rightarrow \mathcal{G}$ est le crochet d'alg\`{e}bre de Lie sur $\mathcal{G}$. On a la d\'{e}finition suivante :
\begin{defi} Le crochet de Schouten alg\'{e}brique est la structure d'alg\`{e}bre de Lie gradu\'{e}e $[, ]^{\mu}$, sur l'alg\`{e}bre ext\'{e}rieure, $\Lambda\mathcal{G} = \bigoplus_{p\geq -1}\Lambda^{p+1}\mathcal{G}$, de $\mathcal{G}$ qui :\\
(i) s'annule si l'un des arguments est dans {\bf K},\\
(ii) \'{e}tend le crochet de Lie $\mu$, i.e
$$[x, y]^{\mu} = \mu(x, y), \forall x, y \in \mathcal{G},$$
(iii)satisfait la r\`{e}gle  suivante sur le degr\'{e} :
$$[X, Y]^{\mu}\in \Lambda^{p+q+1}\mathcal{G},$$
si $X \in \Lambda^{p+1}\mathcal{G}$ et $Y \in \Lambda^{q+1}\mathcal{G}$,\\
(iv) satisfait l'anti-commutativit\'{e} gradu\'{e}e, i.e
$$[X, Y]^{\mu} = - (-1)^{pq}[Y, X]^{\mu},$$
si $X \in \Lambda^{p+1}\mathcal{G}$ et $Y \in \Lambda^{q+1}\mathcal{G}$,\\
(v)satisfait la r\`{e}gle de Leibniz gradu\'{e}e
$$ [X, Y\wedge Z]^{\mu} = [X, Y]^{\mu}\wedge Z + (-1)^{p(q+1)}Y\wedge[X, Z]^{\mu},$$
si $X \in \Lambda^{p+1}\mathcal{G}$, $Y \in \Lambda^{q+1}\mathcal{G}$ et $Z \in \Lambda\mathcal{G}$, et\\
(vi) satisfait l'identit\'{e} de Jacobi gradu\'{e}e, i.e
$$(-1)^{pr}[[X, Y]^{\mu}, Z]^{\mu} + (-1)^{pq}[[Y, Z]^{\mu}, X]^{\mu} + (-1)^{qr}[[Z, X]^{\mu}, Y]^{\mu} = 0,$$
si $X \in \Lambda^{p+1}\mathcal{G}$, $Y \in \Lambda^{q+1}\mathcal{G}$ et $Z \in \Lambda^{r+1}\mathcal{G}$.
\end{defi}
{\bf Remarque 2.1:} On peut d\'{e}finir un tel crochet sur $\Lambda\mathcal{G}$ pour tout \'{e}l\'{e}ment $\mu \in Hom(\Lambda^{2}\mathcal{G}, \mathcal{G})$; l'anti-commutativit\'{e} gradu\'{e}e et la r\`{e}gle de Leibniz gradu\'{e}e restent en vigueur, mais en g\'{e}n\'{e}ral le crochet $[, ]^{\mu}$ ne v\'{e}rifie pas l'identit\'{e} de Jacobi gradu\'{e}e et il la v\'{e}rifie si et seulement si $(\mathcal{G}, \mu)$ est une alg\`{e}bre de Lie (\cite{Kz2}). Dans toute la suite le terme crochet signifiera un \'{e}l\'{e}ment de $Hom(\Lambda^{2}\mathcal{G}, \mathcal{G})$; il sera un crochet de Lie s'il v\'{e}rifie l'identit\'{e} de Jacobi.\\
On a le r\'{e}sultat suivant
\begin{prop} Soit $(\mathcal{G}, \mu)$ une alg\`{e}bre de Lie. Alors, $\forall x \in \mathcal{G}$, $Y\in \Lambda\mathcal{G}$, on a
$$\varepsilon_{[x, Y]^{\mu}} = [ad_{x}, \varepsilon_{Y}],$$
o\`{u} $ad : x\in \mathcal{G} \rightarrow ad_{x}\in End\mathcal{G}$ est la repr\'{e}sentation adjointe de $\mathcal{G}$, d\'{e}finie, pour $y\in \mathcal{G}$, par $ad_{x}y = \mu(x, y)$ et $[, ]$ d\'{e}signe le crochet commutateur des endormorphismes.\\
\end{prop}
La d\'{e}monstration de cette proposition est une traduction de la r\`{e}gle de Leibniz gradu\'{e}e du crochet $[, ]^{\mu}$.  $\bigtriangleup$
\subsection{Cohomologie d'alg\`{e}bre de Lie}
Soit $(\mathcal{G}, \mu)$ une alg\`{e}bre de Lie et soit $M$ un $\mathcal{G}$-module, c'est-\`{a}-dire que $\mathcal{G}$ agit sur $M$. Par exemple $\mathcal{G}$ agit sur elle-m\^{e}me (plus g\'{e}n\'{e}ralement sur son alg\`{e}bre tensorielle) par la repr\'{e}sentation adjointe.\\
$\mathcal{G}$ agit sur son alg\`{e}bre tensorielle de la mani\`{e}re suivante ; pour des \'{e}l\'{e}ments d\'{e}composables, $y_{1}\otimes y_{2}\otimes....\otimes y_{n}\in \bigotimes^{n}\mathcal{G}$,
$$x.(y_{1}\otimes y_{2}\otimes....\otimes y_{n}) = \sum_{i=1}^{n}y_{1}\otimes y_{2}\otimes....\otimes(ad_{x}y_{i})\otimes...\otimes y_{n}.$$
\begin{defi} L'espace vectoriel des applications k-lin\'{e}aires antisym\'{e}triques sur $\mathcal{G}$ \`{a} valeurs dans M est appel\'{e}e l'espace des k-cochaines de $\mathcal{G}$ \`{a} valeurs dans M.
\end{defi}

D\'{e}signons par $\mathcal{C}^{k}(\mathcal{G}, M)$, l'espace vectoriel des k-cochaines de $\mathcal{G}$ \`{a} valeurs dans M et $\mathcal{C}(\mathcal{G}, M) = \bigoplus_{k\geq 0}\mathcal{C}^{k}(\mathcal{G}, M)$.
\begin{defi} L'op\'{e}rateur cobord de Chevalley-Eilenberg de $\mathcal{G}$ \`{a} valeurs dans M, not\'{e} $\delta_{\mu} : \mathcal{C}(\mathcal{G}, M)\rightarrow \mathcal{C}(\mathcal{G}, M)$, est l'application lin\'{e}aire de degr\'{e} $1$ d\'{e}finie par :
$$
\begin{array}{ccc}
(\delta_{\mu}\alpha)(x_{0},...,x_{k}) & = & \sum_{i=0}^{k}(-1)^{i}x_{i}.(\alpha(x_{0},...,\hat{x_{i}},...,x_{k}))\\
&  &\\
& & + \sum_{i<j}(-1)^{i+j}\alpha(\mu(x_{i}, x_{j}), x_{0},...,\hat{x_{i}},...,\hat{x_{j}},...,x_{k})
\end{array}
$$
pour $\alpha \in \mathcal{C}^{k}(\mathcal{G}, M)$, $x_{i}\in \mathcal{G}, i=0,1,...,k$; $x_{i}.m$ d\'{e}signant l'action de $x_{i}\in \mathcal{G}$ sur $m\in M$ et $\hat{x_{i}}$ indiquant l'omission de l'argument $x_{i}$.
\end{defi}
{\bf Remarque 2.2:} On peut d\'{e}finir un tel op\'{e}rateur $\delta_{\mu}$ pour tout \'{e}l\'{e}ment $\mu \in Hom(\Lambda^{2}\mathcal{G}, \mathcal{G})$. En g\'{e}n\'{e}ral $\delta_{\mu}^{2} \neq 0$ et $\delta_{\mu}^{2} = 0$ si et seulement si $(\mathcal{G}, \mu)$ est une alg\`{e}bre de Lie (\cite{Kz1}).
\begin{defi} Une k-cochaine $\alpha$ est appel\'{e}e un k-cocycle si $\delta_{\mu}\alpha = 0$. Une k-cochaine $\alpha$ est appel\'{e}e un k-cobord si il existe une (k-1)-cochaine $\beta$, telle que $\alpha = \delta_{\mu}\beta$.
\end{defi}

Ainsi, comme on le voit, tout k-cobord est un k-cocycle; ce qui nous conduit \`{a} la d\'{e}finition suivante :
\begin{defi} Le quotient de l'espace vectoriel des k-cocycles par l'espace vectoriel des k-cobords est appel\'{e} le $k^{ieme}$ espace de cohomologie de $\mathcal{G}$ \`{a} valeurs dans M et not\'{e} par $H^{k}(\mathcal{G}, M)$.
\end{defi}
{\bf Remarque 2.3:} Les 0-cocycles de $\mathcal{G}$ \`{a} valeurs dans M sont les \'{e}l\'{e}ments invariants dans $M$, i.e les  \'{e}l\'{e}ments $m\in M$ tels que $x.m = 0$, pour tous $x\in \mathcal{G}$.
\subsection{La repr\'{e}sentation co-adjointe}
On introduit \`{a} pr\'{e}sent la d\'{e}finition de la repr\'{e}sentation co-adjointe d'une alg\`{e}bre de Lie sur son espace vectoriel dual.\\
Soit $(\mathcal{G}, \mu)$ une alg\`{e}bre de Lie et soit $\mathcal{G}^{*}$ son espace vectoriel dual. Pour $x\in \mathcal{G}$, posons
$$ad_{x}^{*} = - ^{t}(ad_{x}).$$
Comme on le voit par d\'{e}finition, $ad_{x}^{*}\in End(\mathcal{G}^{*})$ et satisfait la relation
$$<ad_{x}^{*}\xi, x> = - <\xi, ad_{x}y> = - <\xi, \mu(x, y)>,$$
pour $x\in \mathcal{G}$, $\xi\in \mathcal{G}^{*}$. On montre facilement que l'application $x\in \mathcal{G} \rightarrow ad_{x}^{*}\in End(\mathcal{G}^{*})$ est une repr\'{e}sentation de $\mathcal{G}$ dans $\mathcal{G}^{*}$. D'o\`{u}
\begin{defi} La repr\'{e}sentation $x \rightarrow ad_{x}^{*}$ de $\mathcal{G}$ dans $\mathcal{G}^{*}$ est appel\'{e}e la repr\'{e}sentation co-adjointe de $\mathcal{G}$.
\end{defi}

De la proposition 2.1 et de la d\'{e}finition pr\'{e}c\'{e}dente, nous obtenons le r\'{e}sultat suivant :
\begin{prop} Soit $(\mathcal{G}, \mu)$ une alg\`{e}bre de Lie. Alors, $\forall x \in \mathcal{G}$, $Y\in \Lambda\mathcal{G}$, on a
$$i_{[x, Y]^{\mu}} = [ad_{x}^{*}, i_{Y}].$$
\end{prop}
\section{Quasi-big\`{e}bres de Lie}
Les quasi-big\`{e}bres de Lie (\cite{Dr2}) (appel\'{e}es quasi-big\`{e}bres jacobiennes dans (\cite{Ba1}, \cite{BK}, \cite{KS})) sont les limites classiques des alg\`{e}bres quasi-Hopf (\cite{Dr2}), introduites par Drinfeld.
\subsection{D\'{e}finitions et notations}
\begin{defi} Une quasi-big\`{e}bre de Lie est un quadruplet $(\mathcal{G}, \mu, \gamma, \phi)$ o\`{u} $\mathcal{G}$ est un espace vectoriel muni d'un crochet $\mu \in Hom(\Lambda^{2}\mathcal{G}, \mathcal{G})$, d'un co-crochet $\gamma \in Hom(\mathcal{G}, \Lambda^{2}\mathcal{G})$ et d'un \'{e}l\'{e}ment $\phi\in \Lambda^{3}\mathcal{G}$ tels que :\\

3.1. $(\mathcal{G}, \mu)$ est une alg\`{e}bre de Lie;\\

3.2. $\gamma$ est un 1-cocycle de l'alg\`{e}bre de Lie $(\mathcal{G}, \mu)$, \`{a} valeurs dans $\Lambda^{2}\mathcal{G}$ pour l'action adjointe d\'{e}finie par $\mu$, i.e  $\delta_{\mu}\gamma = 0;$\\

3.3. $\frac{1}{2}Alt(\gamma\otimes 1)\gamma(x) = (\delta_{\mu}\phi)(x),  \forall x \in \mathcal{G}$;\\

3.4. $Alt(\gamma\otimes 1\otimes 1)(\phi) = 0;$\\

o\`{u} Alt est l'op\'{e}rateur alternateur d\'{e}fini sur l'alg\`{e}bre tensorielle de $\mathcal{G}$ par
$$Alt(X_{1}\otimes....\otimes x_{n}) = \sum_{\sigma}sign(\sigma)x_{\sigma(1)}\otimes... \otimes x_{\sigma(n)},$$
$x_{i}\in \mathcal{G}, i = 1,...,n$, $\sigma$ \'{e}tant une permutation de $\{1,...,n\}$ et $sign(\sigma)$ la signature de la permutation $\sigma$.
\end{defi}
{\bf Remarque 3.1 :} \\
1- Dans le cas o\`{u} $\phi = 0$, le triplet $(\mathcal{G}, \mu, \gamma)$ satisfaisant les conditions ci-dessus, n'est rien d'autre qu'une big\`{e}bre de Lie (\cite{Dr1}).\\
2- La condition 2.3 signifie que $\gamma$ ne v\'{e}rifie pas l'identit\'{e} co-Jacobi et donc son transpos\'{e} n'est pas un crochet de Lie sur $\mathcal{G}^{*}$. Ainsi, contrairement \`{a} la notion de big\`{e}bre de Lie, la notion de quasi-big\`{e}bre de Lie n'est pas auto-duale, l'objet dual est appel\'{e} une big\`{e}bre quasi-Lie (\cite{Dr2}) ou quasi-big\`{e}bre co-jacobienne (\cite{Ba1}, \cite{KS}). Dans ce travail nous consid\'{e}rerons pour des fins d'usage, l'oppos\'{e} du transpos\'{e} de $\gamma$ comme \'{e}tant le crochet sur $\mathcal{G}^{*}$ et nous le noterons aussi par $\gamma$ pour simplicit\'{e} d'\'{e}criture, i.e
$$<\gamma(x), \xi\wedge \eta> = - <x, \gamma(\xi, \eta)>, \forall x\in \mathcal{G}, \forall \xi, \eta\in \mathcal{G}^{*}.$$
3- Dans la condition 2.3, $\phi : {\bf K}\rightarrow \Lambda^{3}\mathcal{G}$ est consid\'{e}r\'{e} comme une $0$-forme sur $\mathcal{G}$ \`{a} valeurs dans $\Lambda^{3}\mathcal{G}$, tandis que si nous consid\'{e}rons $\phi : \Lambda^{3}\mathcal{G}^{*}\rightarrow {\bf K}$ comme une 3-forme sur $\mathcal{G}^{*}$ \`{a} valeurs dans {\bf K}, alors la condition 2.4 de la d\'{e}finition ci-dessus est \'{e}quivalente \`{a} $\delta_{\gamma}\phi = 0$.\\
4- La condition 2.2 s'\'{e}crit explicitement sous la forme
$$\gamma(\mu(x, y)) = (ad_{x}^{\mu}\otimes 1 + 1\otimes ad_{x}^{\mu})\gamma(y) - (ad_{y}^{\mu}\otimes 1 + 1\otimes ad_{y}^{\mu})\gamma(x),$$
o\`{u} $ad_{x}^{\mu}y = \mu(x, y), \forall x, y \in \mathcal{G}$, ou de fa\c{c}on \'{e}quivalente
$$\mu(\gamma(\xi, \eta)) = (ad_{\xi}^{\gamma}\otimes 1 + 1\otimes ad_{\xi}^{\gamma})\mu(\eta) - (ad_{\eta}^{\gamma}\otimes 1 + 1\otimes ad_{\eta}^{\gamma})\mu(\xi),$$
o\`{u} $ad_{\xi}^{\gamma}\eta = \gamma(\xi, \eta), \forall \xi, \eta \in \mathcal{G}^{*}$, et l'application $\mu$ et sa transpos\'{e}e sont not\'{e}es par $\mu$ pour simplicit\'{e} d'\'{e}criture.

La donn\'{e}e d'une structure de quasi-big\`{e}bre de Lie sur $\mathcal{G}$ d\'{e}termine une unique structure d'alg\`{e}bre de Lie $[ , ]_{\mathcal{D}}$ sur l'espace vectoriel $\mathcal{D} = \mathcal{G} \oplus \mathcal{G}^{*}$ qui laisse invariant le produit scalaire canonique sur $\mathcal{D}$,
$$ <\xi + x, y + \eta> = <\xi, y> + <x, \eta>, \hspace{0.5cm}\forall \xi + x \in \mathcal{D}^{*},  \forall y + \eta \in \mathcal{D},$$
en posant
$$[x, y]_{\mathcal{D}} = \mu(x, y),$$
$$[x, \xi]_{\mathcal{D}} = - ad_{\xi}^{\gamma *}x + ad_{x}^{\mu *}\xi,$$
$$[\xi, \eta]_{\mathcal{D}} = \phi(\xi, \eta) + \gamma(\xi, \eta),$$
o\`{u} $<ad_{x}^{\mu *}\xi, y> = - <\xi, \mu(x, y)>$, $<ad_{\xi}^{\gamma *}x, \eta> = - <x, \gamma(\xi, \eta)>$ et $\phi(\xi, \eta) = i_{\xi\wedge\eta}\phi$, pour tous $x, y \in \mathcal{G}, \xi, \eta \in \mathcal{G}^{*}$.\\
Comme on le voit, $(\mathcal{G}, \mu)$ est une sous-alg\`{e}bre de Lie isotrope de $(\mathcal{D}, [ , ]_{\mathcal{D}})$, alors que $\mathcal{G}^{*}$ ne l'est pas, il est juste un sous-espace isotrope de $\mathcal{D}$, \`{a} cause de l'existence de $\phi$. En g\'{e}n\'{e}ral, on montre dans (\cite{KS}) que les structures de quasi-big\`{e}bre de Lie sur $\mathcal{G}$ sont en correspondance biunivoque avec les structures d'alg\`{e}bre de Lie sur $\mathcal{D} = \mathcal{G} \oplus \mathcal{G}^{*}$ laissant invariant le produit scalaire canonique, dont $\mathcal{G}$ est une sous-alg\`{e}bre de Lie. Dans ces conditions le couple $(\mathcal{D}, \mathcal{G})$, avec le produit scalaire canonique sur $\mathcal{D}$, est appel\'{e} un {\bf couple de Manin} (\cite{Dr2}). Plus pr\'{e}cisement
\begin{defi} Un couple de Manin consiste en un couple $(\mathcal{D}, \mathcal{G})$, o\`{u} $\mathcal{D}$ est une alg\`{e}bre de Lie munie d'un produit scalaire invariant non d\'{e}g\'{e}n\'{e}r\'{e} et $\mathcal{G}$ une sous-alg\`{e}bre de Lie isotrope de dimension maximale de $\mathcal{D}$.
\end{defi}

L'\'{e}tude des quasi-big\`{e}bres de Lie est rendue facile gr\^{a}ce au twisting (\cite{Dr2}, \cite{KS}) appel\'{e} modification dans (\cite{Ba1}), qui consiste \`{a} construire de nouvelles structures de quasi-big\`{e}bre de Lie sur $\mathcal{G}$ \`{a} partir d'une d\'{e}j\`{a} connue; ce qui permet de les \'{e}tudier en termes de classes d'\'{e}quivalence (\cite{KS}), en montrant que les classes d'\'{e}quivalence modulo twisting sont en correspondance biunivoque avec les couples de Manin.
\begin{defi} Soit $(\mathcal{G}, \mu, \gamma, \phi)$ une quasi-big\`{e}bre de Lie. $\mathcal{D} = \mathcal{G} \oplus \mathcal{G}^{*}$ muni du crochet de Lie $[ , ]_{\mathcal{D}}$ d\'{e}fini ci-dessus est appel\'{e} le double de la quasi-big\`{e}bre de Lie donn\'{e}e, et not\'{e} $\mathcal{G} \bowtie \mathcal{G}^{*}$.
\end{defi}
\begin{defi} Une quasi-big\`{e}bre de Lie $(\mathcal{G}, \mu, \gamma, \phi)$ est dite exacte ou cobord si il existe un \'{e}l\'{e}ment ${\bf r} \in \Lambda^{2}\mathcal{G}$ tel que le $1$-cocycle $\gamma : \mathcal{G}\rightarrow \Lambda^{2}\mathcal{G}$ soit le cobord de {\bf r}, i.e
$$\gamma(x) = (\delta_{\mu}{\bf r})(x) = [x, {\bf r}]^{\mu} = - [{\bf r}, x]^{\mu}, \forall x\in \mathcal{G},$$
et
$$\phi = - \frac{1}{2}[{\bf r}, {\bf r}]^{\mu}.$$
\end{defi}

Montrons \`{a} pr\'{e}sent que le double de toute quasi-big\`{e}bre de Lie est muni, en plus de la structure d'alg\`{e}bre de Lie d\'{e}finie par $[ , ]_{\mathcal{D}}$, d'une structure canonique de quasi-big\`{e}bre de Lie exacte (\cite{BK}).
\begin{thm} Soit $\mathcal{D} = \mathcal{G} \bowtie \mathcal{G}^{*}$ le double d'une quasi-big\`{e}bre de Lie $(\mathcal{G}, \mu, \gamma, \phi)$. Soit $(e_{i})$ une base de $\mathcal{G}$ et $(\xi^{i})$ la base duale de $\mathcal{G}^{*}$. Posons
$$ {\bf r} = \frac{1}{2}\sum_{i}e_{i}\wedge \xi^{i}.$$
Alors $(\mathcal{D}, {\bf r})$ est une quasi-big\`{e}bre de Lie exacte et est appel\'{e}e la quasi-big\`{e}bre de Lie double de $(\mathcal{G}, \mu, \gamma, \phi)$.
\end{thm}
Nous avons le r\'{e}sultat suivant \cite{Ba1}:
\begin{prop} Soit $(\mathcal{G}, \mu, \gamma, \phi)$ une quasi-big\`{e}bre de Lie. Alors nous avons les relations suivantes :\\

3.5. $ad_{\mu(x, y)}^{\mu *} = [ad_{x}^{\mu *}, ad_{y}^{\mu *}], \forall x, y \in \mathcal{G}$;\\

3.6. $ad_{\xi}^{\gamma *}\mu(x, y) = \mu(ad_{\xi}^{\gamma *}x, y) + \mu(x, ad_{\xi}^{\gamma *}y) + ad_{ad_{y}^{\mu *}\xi}^{\gamma *}x - ad_{ad_{x}^{\mu *}\xi}^{\gamma *}y \\
\forall x, y \in \mathcal{G}, \forall \xi\in \mathcal{G}^{*};$\\

3.7. $ad_{\gamma(\xi, \eta)}^{\gamma *}x = [ad_{\xi}^{\gamma *}, ad_{\eta}^{\gamma *}](x) + ad_{x}^{\mu}\phi(\xi, \eta) - \phi(ad_{x}^{\mu *}\xi, \eta) - \phi(\xi, ad_{x}^{\mu *}\eta) ,\\
 \forall x \in \mathcal{G}, \forall \xi, \eta \in \mathcal{G}^{*};$\\

3.8. $ad_{x}^{\mu *}\gamma(\xi, \eta) = \gamma(ad_{x}^{\mu *}\xi, \eta) + \gamma(\xi, ad_{x}^{\mu *}\eta) + ad_{ad_{\eta}^{\gamma *}x}^{\mu *}\xi - ad_{ad_{\xi}^{\gamma *}x}^{\mu *}\eta,\\
 \forall x\in \mathcal{G}, \forall \xi, \eta\in \mathcal{G}^{*};$\\

3.9. $\oint\gamma(\gamma(\xi, \eta), \zeta) = - \oint ad_{\phi(\xi, \eta)}^{\mu *}\zeta, \forall \xi, \eta, \zeta \in \mathcal{G}^{*};$\\

3.10. $\oint\phi(\gamma(\xi, \eta), \zeta) = \oint ad_{\xi}^{\gamma *}\phi(\xi, \eta), \forall \xi, \eta, \zeta \in \mathcal{G}^{*}.$\\
o\`{u} $[, ]$ d\'{e}signe le crochet commutateur des endormorphismes et $\oint$ d\'{e}signe la somme sur les permutations circulaires des \'{e}l\'{e}ments $\xi, \eta, \zeta \in \mathcal{G}^{*}.$
\end{prop}
{\bf D\'{e}monstration :} La preuve de ces diff\'{e}rentes relations est une cons\'{e}quence directe de l'identit\'{e} de Jacobi pour le crochet d'alg\`{e}bre de Lie $[ , ]_{\mathcal{D}}$ d\'{e}fini sur $\mathcal{D} = \mathcal{G} \bowtie \mathcal{G}^{*}$. $\bigtriangleup$

Nous avons les \'{e}quivalences suivantes :\\
- la relation (3.5) traduit le fait que $\mu$ d\'{e}fini un crochet de Lie sur $\mathcal{G}$;\\
- les relations (3.6) et (3.8) sont \'{e}quivalentes \`{a} la condition de 1-cocycle pour $\gamma$, $\delta_{\mu}\gamma = 0$;\\
- les relations (3.7) et (3.9) sont \'{e}quivalentes \`{a} la condition $\frac{1}{2}Alt(\gamma\otimes 1)\gamma(x) = (\delta_{\mu}\phi)(x)$;\\
- la relation (3.10) est \'{e}quivalente \`{a} la condition (3.4), i.e $Alt(\gamma\otimes 1\otimes 1)(\phi) = 0.$

Les relations (3.6) et (3.7) s'\'{e}tendent aisement sur $\Lambda\mathcal{G}$ gr\^{a}ce au r\'{e}sultat suivant :
\begin{prop} Pour tous $X = x_{1}\wedge...\wedge x_{m}$ et $Y = y_{1}\wedge...\wedge y_{n}$ dans $\Lambda\mathcal{G}$, et pour tous $\xi, \eta \in \mathcal{G}^{*}$, on a
$$\begin{array}{ccc}
 ad_{\xi}^{\gamma *}[X, Y]^{\mu} & = & [ad_{\xi}^{\gamma *}X, Y]^{\mu} + [X, ad_{\xi}^{\gamma *}Y]^{\mu}\\
  &   &\\
  &   & + \sum_{j=1}^{n}(-1)^{j+1}ad_{ad_{{y_{j}}}^{\mu *}\xi}^{\gamma *}X\wedge \hat{Y_{j}} \\
  &   & \\
  &    & + (-1)^{|X|} \sum_{i=1}^{m}(-1)^{i+1}\hat{X_{i}}\wedge ad_{ad_{{x_{i}}}^{\mu *}\xi}^{\gamma *}Y,
\end{array}
$$
et
$$\begin{array}{ccc}
ad_{\gamma(\xi, \eta)}^{\gamma *}X & = & [ad_{\xi}^{\gamma *}, ad_{\eta}^{\gamma *}](X) -  ad_{\phi(\xi, \eta)}^{\mu}X\\
 &   &  \\
 &    & + \sum_{i=1}^{m}(-1)^{i}(\phi(ad_{x_{i}}^{\mu *}\xi, \eta) + \phi(\xi, ad_{x_{i}}^{\mu *}\eta))\wedge\hat{X_{i}},
 \end{array}
$$
o\`{u} $\hat{X}_{i} = x_{1}\wedge...\wedge \hat{x_{i}}\wedge...\wedge x_{m}$ pour $1\leq i\leq m$ et de mani\`{e}re similaire pour $\hat{Y}_{j}$, pour $1\leq j\leq n$.
\end{prop}
{\bf D\'{e}monstration :} Ces deux relations se d\'{e}montrent facilement par r\'{e}currence sur les degr\'{e}s de $X$ et $Y$, en utilisant la r\`{e}gle de Leibniz gradu\'{e}e du crochet $[, ]^{\mu}$ et la propri\'{e}t\'{e} de d\'{e}rivation de $ad_{\xi}^{\gamma *}$ sur $(\Lambda\mathcal{G}, \wedge)$ pour tout $\xi \in \mathcal{G}^{*}$. La premi\`{e}re relation est \'{e}nonc\'{e}e dans (\cite{Lu}) dans le cas des big\`{e}bres de Lie. ${\bf \bigtriangleup}$
\subsection{Exemples}
\begin{ex} Toute big\`{e}bre de Lie est une quasi-big\`{e}bre de Lie; il suffit de prendre $\phi = 0$.
\end{ex}
\begin{ex} Une large classe d'exemples de quasi-big\`{e}bre de Lie est fournie par les quasi-big\`{e}bres de Lie exactes, il suffit de choisir ${\bf r} \in \Lambda^{2}\mathcal{G}$.
\end{ex}
\begin{ex} Soit $(\mathcal{G}, \mu)$ une alg\`{e}bre de Lie. Alors tout \'{e}l\'{e}ment ${\bf r} \in \mathcal{G}\otimes \mathcal{G}$ de partie antisym\'{e}trique $ad^{\mu}$-invariante, d\'{e}finie une structure de quasi-big\`{e}bre de Lie en posant
$$\gamma = \delta_{\mu}{\bf a}, \hspace{1cm}\phi = - \frac{1}{2}([{\bf a}, {\bf a}]^{\mu} + [{\bf s}, {\bf s}]^{\mu}),$$
o\`{u} {\bf a} (resp. {\bf s}) est la partie antisym\'{e}trique (resp. sym\'{e}trique) de {\bf r}. Une telle structure est dite quasitriangulaire (\cite{Ba1}, \cite{BK}).
\end{ex}
\begin{ex} Soit $(\mathcal{D}, \mathcal{G})$ un couple de Manin; alors tout choix d'un sous-espace suppl\'{e}mentaire
isotrope de $\mathcal{G}$ dans $\mathcal{D}$ d\'{e}finit une structure de quasi-big\`{e}bre de Lie sur $\mathcal{G}$.
\end{ex}
\subsection{Les op\'{e}rateurs de cohomologie sur une quasi-big\`{e}bre de Lie}
Soit $(\mathcal{G}, \mu, \gamma, \phi)$ une quasi-big\`{e}bre de Lie. Alors le crochet $\mu$ permet de d\'{e}finir sur $\Lambda\mathcal{G}^{*}$ l'op\'{e}rateur cobord de  Chevalley-Eilenberg (\`{a} coefficients triviaux) et son transpos\'{e} sur $\Lambda\mathcal{G}$; notons-les respectivement par $d_{\mu}$ et $\partial_{\mu}$:
$$ d_{\mu} : \Lambda^{k}\mathcal{G}^{*} \rightarrow \Lambda^{k+1}\mathcal{G}^{*}$$
$$ \partial_{\mu} : \Lambda^{k}\mathcal{G} \rightarrow \Lambda^{k-1}\mathcal{G}$$
o\`{u}
$$(d_{\mu}\xi)(x_{1}\wedge...\wedge x_{k+1}) = \sum_{i<j}(-1)^{i+j}\xi(\mu(x_{i}, x_{j})\wedge x_{1}\wedge...\wedge\hat{x_{i}}\wedge...\wedge\hat{x_{j}}\wedge...\wedge x_{k+1}),$$
$$\partial_{\mu}(x_{1}\wedge...\wedge x_{k+1}) = \sum_{i<j}(-1)^{i+j}\mu(x_{i}, x_{j})\wedge x_{1}\wedge...\wedge\hat{x_{i}}\wedge...\wedge\hat{x_{j}}\wedge...\wedge x_{k+1},$$
pour $\xi \in \Lambda^{k}\mathcal{G}^{*}$, $x_{i} \in \mathcal{G}, i = 1,...,k+1$.\\
De mani\`{e}re similaire, d\'{e}finissons les op\'{e}rateurs $d_{\gamma}$ et $\partial_{\gamma}$ associ\'{e}s au crochet $\gamma$; ils sont d\'{e}finis respectivement sur $\Lambda\mathcal{G}$ et $\Lambda\mathcal{G}^{*}$ comme suit :
$$ d_{\gamma} : \Lambda^{k}\mathcal{G} \rightarrow \Lambda^{k+1}\mathcal{G}$$
$$ \partial_{\gamma} : \Lambda^{k}\mathcal{G}^{*} \rightarrow \Lambda^{k-1}\mathcal{G}^{*}$$
o\`{u}
$$(d_{\gamma}X)(\xi_{1}\wedge...\wedge \xi_{k+1}) = \sum_{i<j}(-1)^{i+j}X(\gamma(\xi_{i}, \xi_{j})\wedge \xi_{1}\wedge...\wedge\hat{\xi_{i}}\wedge...\wedge\hat{\xi_{j}}\wedge...\wedge \xi_{k+1}),$$
$$\partial_{\gamma}(\xi_{1}\wedge...\wedge \xi_{k+1}) = \sum_{i<j}(-1)^{i+j}\gamma(\xi_{i},\xi_{j})\wedge \xi_{1}\wedge...\wedge\hat{\xi_{i}}\wedge...\wedge\hat{\xi_{j}}\wedge...\wedge \xi_{k+1},$$
pour $X \in \Lambda^{k}\mathcal{G}$, $\xi_{i} \in \mathcal{G}^{*}, i = 1,...,k+1$.\\
{\bf Remarque 3.2 :} Pour $k = 1$, $d_{\gamma}$ n'est rien d'autre que le 1-cocycle $\gamma$.\\
Il est connu  que l'op\'{e}rateur $d_{\mu}$ est de carr\'{e} nul, mais  $d_{\gamma}^{2}\neq 0$ du fait que $\gamma$ ne d\'{e}finit pas un crochet de Lie sur $\mathcal{G}^{*}$. Par ailleurs, les op\'{e}rateurs $d_{\mu}$ et $\partial_{\mu}$ satisfont les propri\'{e}t\'{e}s suivantes (\cite{Kz2}, \cite{Lu}):
$$d_{\mu}(A\wedge B) = (d_{\mu}A)\wedge B + (-1)^{|A|}A\wedge d_{\mu}B),$$
$$\partial_{\mu}(X\wedge Y) = (\partial_{\mu}X)\wedge Y + (-1)^{|X|}X\wedge \partial_{\mu}Y) + (-1)^{|X|}[X, Y]^{\mu},$$
$$\partial_{\mu}[X, Y]^{\mu} = [\partial_{\mu}X, Y]^{\mu} + (-1)^{(|X|-1)}[X, \partial_{\mu}Y]^{\mu},$$
pour tous $X, Y \in \Lambda\mathcal{G}$ et $A, B \in \Lambda\mathcal{G}^{*}$. Les op\'{e}rateurs $d_{\gamma}$ et $\partial_{\gamma}$ satisfont les relations similaires.\\
Les lemmes suivants seront d'une grande utilit\'{e} dans la suite :
\begin{lem} Pour une quasi-big\`{e}bre de Lie donn\'{e}e $(\mathcal{G}, \mu, \gamma, \phi)$, on a :\\

3.11. $ad_{x}^{\mu *} = [d_{\mu}, i_{x}], \forall x \in  \mathcal{G};$\\

3.12. $[d_{\mu}, ad_{x}^{\mu *}] = 0, \forall x \in  \mathcal{G}$, i.e $d_{\mu}(ad_{x}^{\mu *}\xi) = ad_{x}^{\mu *}(d_{\mu}\xi), \forall x \in  \mathcal{G}, \forall \xi\in \mathcal{G}^{*};$\\

3.13. $i_{d_{\mu}(ad_{x}^{\mu *}\xi)} = [ad_{x}^{\mu}, i_{d_{\mu}\xi}], \forall x \in  \mathcal{G}, \forall \xi\in \mathcal{G}^{*}$;\\

3.14. $d_{\gamma}(ad_{\xi}^{\gamma *}x) = ad_{\xi}^{\gamma *}(d_{\gamma}x) + ad_{x}^{\mu}(i_{\xi}\phi) - i_{ad_{x}^{\mu *}\xi}\phi, \forall x \in  \mathcal{G}, \forall \xi\in \mathcal{G}^{*};$\\

3.15. $i_{\gamma(\xi, \eta)}\phi = ad_{\xi}^{\gamma *}i_{\eta}\phi - ad_{\eta}^{\gamma *}i_{\xi}\phi + d_{\gamma}(\phi(\xi, \eta)), \forall \xi, \eta\in \mathcal{G}^{*}$.
\end{lem}
{\bf D\'{e}monstration :}\\
(3.11) Par d\'{e}finition, comme $\partial_{\mu}(x) =0, \forall x \in  \mathcal{G}$ pour cause de degr\'{e}, on a:\\
$$ad_{x}^{\mu}Y = [x, Y]^{\mu} = - \partial_{\mu}(x\wedge Y) - x\wedge\partial_{\mu}Y = - [\partial_{\mu}, \varepsilon_{x}](Y), \forall Y \in  \Lambda\mathcal{G}$$
i.e
$$ad_{x}^{\mu} = - [\partial_{\mu}, \varepsilon_{x}];$$
d'o\`{u} par transpositition $ad_{x}^{\mu *} = [d_{\mu}, i_{x}], \forall x \in  \mathcal{G};$\\
(3.12) et (3.13) sont des cons\'{e}quences directes de (3.11);\\
(3.14) est une cons\'{e}quence directe de la relation (3.9);\\
(3.15) est une cons\'{e}quence directe de la relation (3.10).\\
Ce qui ach\`{e}ve la d\'{e}monstration du lemme. $\bigtriangleup$

\begin{lem} Pour une quasi-big\`{e}bre de Lie donn\'{e}e $(\mathcal{G}, \mu, \gamma, \phi)$, on a pour tout $\forall x \in  \mathcal{G}$, pour tous $\xi, \eta \in \mathcal{G}^{*}$ et pour tout $Y = y_{1}\wedge...\wedge y_{m}\in\Lambda\mathcal{G}$:\\
$$\begin{array}{ccc}
\sum_{i=1}^{m}(-1)^{i+1}ad_{ad_{{y_{i}}}^{\mu *}\xi}^{\gamma *}x\wedge \hat{Y_{i}} & = &
i_{d_{\mu}(\xi)}(d_{\gamma}(x)\wedge Y) - (i_{d_{\mu}(\xi)}(d_{\gamma}(x)))Y \\
&  & \\
&  & - (d_{\gamma}(x))\wedge i_{d_{\mu}(\xi)}Y;\\
&  &\\
\sum_{i=1}^{m}(-1)^{i}(\phi(ad_{y_{i}}^{\mu *}\xi, \eta)\wedge \hat{Y_{i}} & = & i_{d_{\mu}(\xi)}((i_{\eta}\phi)\wedge Y) - (i_{d_{\mu}(\xi)}(i_{\eta}\phi))Y \\
&   &\\
&   &  - (i_{\eta}\phi)\wedge i_{d_{\mu}(\xi)}Y,
\end{array}
$$
o\`{u} $\hat{Y_{i}} = y_{1}\wedge...\wedge \hat{y_{i}}\wedge...\wedge y_{m}$ pour $1\leq i\leq m$.
\end{lem}
La d\'{e}monstration du lemme rel\`{e}ve d'un simple calcul. $\bigtriangleup$

Dans (\cite{Ba2}), on a le r\'{e}sultat suivant :
\begin{prop} Soit $(\mathcal{G}, \mu, \gamma, \phi)$ une quasi-big\`{e}bre de Lie. Alors les op\'{e}rateurs  $d_{\mu}$, $\partial_{\mu}$ et $\partial_{\gamma}$ satisfont les propri\'{e}t\'{e}s suivantes :\\

3.16. $\partial_{\gamma}^{2} + d_{\mu}i_{\phi} + i_{\phi}d_{\mu} - i_{\partial_{\mu}\phi} = 0$;\\

3.17. $\partial_{\gamma}i_{\phi} + i_{\phi}\partial_{\gamma} = 0$;\\

3.18.  $\partial_{\gamma}i_{x} + i_{x}\partial_{\gamma} = - i_{\gamma(x)}$, $\forall x \in \mathcal{G}$.
\end{prop}
On remarque bien que si $\phi = 0$, alors $\partial_{\gamma}^{2} =0$ et par transposition $d_{\gamma}^{2} = 0$.
\begin{defi} L'op\'{e}rateur
$$L =  \partial_{\mu}d_{\gamma} + d_{\gamma}\partial_{\mu} : \Lambda^{k}\mathcal{G} \rightarrow \Lambda^{k}\mathcal{G}$$
est appel\'{e} le {\bf laplacien} de la quasi-big\`{e}bre de Lie $(\mathcal{G}, \mu, \gamma, \phi)$.
\end{defi}
Nous avons le r\'{e}sultat suivant :
\begin{prop} Le laplacien d'une quasi-big\`{e}bre de Lie $(\mathcal{G}, \mu, \gamma, \phi)$ satisfait les propri\'{e}t\'{e}s suivantes :\\

3.19. $L =  \partial_{\mu}d_{\gamma} + d_{\gamma}\partial_{\mu}$  est une d\'{e}rivation de $(\Lambda\mathcal{G}, \wedge)$ de degr\'{e} $0$ et son transpos\'{e} $L^{*} = d_{\mu}\partial_{\gamma} + \partial_{\gamma}d_{\mu}$ est une d\'{e}rivation de $(\Lambda\mathcal{G}^{*}, \wedge)$ de degr\'{e} $0$;\\

3.20. Les op\'{e}rateurs $L$ et $\partial_{\mu}$ commutent, i.e $[L, \partial_{\mu}] = 0$, ou de mani\`{e}re \'{e}quivalente $d_{\mu}$ et $L^{*} = d_{\mu}\partial_{\gamma} + \partial_{\gamma}d_{\mu}$ commutent, i.e $[L^{*}, d_{\mu}] = 0$.
\end{prop}

Dans le cas des big\`{e}bres de Lie, on g\'{e}n\'{e}ralise dans (\cite{Lu}) la condition de 1-cocycle  en montrant que les  op\'{e}rateurs $d_{\mu}$ et $d_{\gamma}$ sont des d\'{e}rivations respectivement de $(\Lambda\mathcal{G}^{*}, [, ]^{\gamma})$ et $(\Lambda\mathcal{G}, [, ]^{\mu})$; le r\'{e}sultat reste vrai dans le cas des quasi-big\`{e}bres de Lie et s'\'{e}nonce comme suit :
\begin{prop} Soit $(\mathcal{G}, \mu, \gamma, \phi)$ une quasi-big\`{e}bre de Lie. Alors, pour $A,B \in\Lambda\mathcal{G}^{*}$ et $X, Y \in \Lambda\mathcal{G}$, on a :
$$d_{\mu}([A, B]^{\gamma}) = [d_{\mu}A, B]^{\gamma} + (-1)^{|A|-1}[A, d_{\mu}B]^{\gamma}$$
$$d_{\gamma}([X, Y]^{\mu}) = [d_{\gamma}X, Y]^{\mu} + (-1)^{|X|-1}[X, d_{\gamma}Y]^{\mu}$$
\end{prop}
{\bf D\'{e}monstration :} Dans le cas o\`{u} $|X| = |Y| = 1$ et $|A| = |B| = 1$, les deux identit\'{e}s se r\'{e}duisent \`{a} la condition de 1-cocycle dans la d\'{e}finition d'une quasi-big\`{e}bre de Lie. Le cas g\'{e}n\'{e}ral se d\'{e}montre par r\'{e}currence sur les degr\'{e}s de $X, Y, A$ et $B$. $\bigtriangleup$

Soit  $(\mathcal{G}, \mu, \gamma, \phi)$ une quasi-big\`{e}bre de Lie telle que $\partial_{\mu}(\phi) \in Im\gamma$, i.e qu'il existe $x_{0} \in \mathcal{G}$ tel que $\gamma(x_{0}) = \partial_{\mu}(\phi)$. Des exemples de telles structures sont fournies par les quasi-big\`{e}bres de Lie exactes o\`{u}
$$\gamma(x) = (\delta_{\mu}{\bf r})(x) = [x, {\bf r}]^{\mu} = - [{\bf r}, x]^{\mu}, \forall x\in \mathcal{G},$$
et
$$\phi = - \frac{1}{2}[{\bf r}, {\bf r}]^{\mu},$$
pour un certain \'{e}l\'{e}ment ${\bf r}\in \Lambda^{2}\mathcal{G}$. En effet, comme $\partial_{\mu}$ est une d\'{e}rivation de $(\Lambda\mathcal{G}, [, ]^{\mu})$, on a
$$\begin{array}{ccc}
\partial_{\mu}(\phi)& = & - \frac{1}{2}\partial_{\mu}([{\bf r}, {\bf r}]^{\mu}) = - \frac{1}{2}([\partial_{\mu}{\bf r}, {\bf r}]^{\mu} - [{\bf r}, \partial_{\mu}{\bf r}]^{\mu})\\
&   &  \\
&  = & - [\partial_{\mu}{\bf r}, {\bf r}]^{\mu} = - \gamma(\partial_{\mu}{\bf r});
\end{array}
$$
d'o\`{u} $\partial_{\mu}(\phi) \in Im\gamma$.\\
Dans \cite{Ba2}, on montre qu'une telle structure de quasi-big\`{e}bre de Lie permet de d\'{e}finir une structure d'alg\`{e}bre quasi-Batalin-Vilkovisky diff\'{e}rentielle (\cite{Ba2}, \cite{G}) sur $\Lambda\mathcal{G}^{*}$ en posant
$$\Delta = \partial_{\gamma} + i_{x_{0}}, \hspace{0.5cm} \delta = d_{\mu}, \hspace{0.5cm} \Phi = i_{\phi}.$$
Plus pr\'{e}cisement on a (\cite{Ba2})
\begin{thm} Les structures de quasi-big\`{e}bre de Lie $(\mathcal{G}, \mu, \gamma, \phi)$ telles que $\partial_{\mu}(\phi) \in Im\gamma$ sont en correspondance bijective avec les structures d'alg\`{e}bre quasi-Batalin-Vilkovisky diff\'{e}rentielle sur $\Lambda\mathcal{G}^{*}$.
\end{thm}

Soit $\xi^{\mu}\in \mathcal{G}^{*}$ et $x^{\gamma}\in \mathcal{G}$ d\'{e}finis respectivement par
$$<\xi^{\mu}, x> = tr(ad_{x}^{\mu}\in End(\mathcal{G})), \hspace{1cm} \forall x \in  \mathcal{G}$$
et
$$<x^{\gamma}, \xi> = tr(ad_{\xi}^{\gamma}\in End(\mathcal{G}^{*})), \hspace{1cm} \forall \xi \in  \mathcal{G}^{*}.$$
On a le r\'{e}sultat suivant :
\begin{lem} Soit $(\mathcal{G}, \mu, \gamma, \phi)$ une quasi-big\`{e}bre de Lie. Alors  les \'{e}l\'{e}ments $\xi^{\mu}$ et $x^{\gamma}$ satisfont les propri\'{e}t\'{e}s suivantes :\\

3.21. $ ad_{x}^{\mu *}\xi^{\mu} = 0, \forall x \in  \mathcal{G}$, ou de mani\`{e}re \'{e}quivalente $d_{\mu}(\xi^{\mu}) =0;$\\

3.22. $<x^{\gamma}, \gamma(\xi, \eta)> = <\xi^{\mu}, \phi(\xi, \eta)> + 2(i_{d_{\mu}(\xi)}i_{\eta}\phi) - 2(i_{d_{\mu}(\eta)}i_{\xi}\phi), \forall \xi, \eta\in \mathcal{G}^{*}$, ou de mani\`{e}re \'{e}quivalente $d_{\gamma}(x^{\gamma}) = -i_{\xi^{\mu}}\phi - 2\partial_{\mu}\phi$;\\

3.23. $<x^{\gamma}, ad_{x}^{\mu *}\xi> = - <\xi^{\mu}, ad_{\xi}^{\gamma *}x> + 2(i_{d_{\mu}(\xi)}d_{\gamma}(x)), \forall x \in  \mathcal{G}, \forall \xi\in \mathcal{G}^{*}.$
\end{lem}
{\bf D\'{e}monstration :} Elle utilise essentiellement les d\'{e}finitions des diff\'{e}rents op\'{e}rateurs et les axiomes d\'{e}finissant la structure de quasi-big\`{e}bre de Lie. La relation (3.21) est une \'{e}vidence; la relation (3.22) suit de (3.7). La relation (3.23) est une cons\'{e}quence de (3.8). $\bigtriangleup$

Dans la th\'{e}orie des big\`{e}bres de Lie, $\xi^{\mu}$ est appel\'{e} le {\bf caract\`{e}re adjoint} de $\mathcal{G}$ et $x^{\gamma}$ le {\bf caract\`{e}re adjoint} de $\mathcal{G}^{*}$ (\cite{Lu}). On montre dans ce cas que le laplacien de la big\`{e}bre de Lie $(\mathcal{G}, \mu, \gamma)$ est
$$L = \frac{1}{2}(ad_{x^{\gamma}}^{\mu} - ad_{\xi^{\mu}}^{\gamma *})$$
et son application duale est
$$L^{*} = \frac{1}{2}(- ad_{x^{\gamma}}^{\mu *} + ad_{\xi^{\mu}}^{\gamma}).$$
\section{Repr\'{e}sentation de $\mathcal{D} = \mathcal{G} \bowtie \mathcal{G}^{*}$ sur $\Lambda\mathcal{G}$}
Dans cette section, on montre qu'il existe une structure de $\mathcal{D}$-module, ou de mani\`{e}re \'{e}quivalente une repr\'{e}sentation de  $\mathcal{D}$, sur $\Lambda\mathcal{G}$, puis on montre que $\Lambda\mathcal{D}$  et $End(\Lambda\mathcal{G})$, comme  $\mathcal{D}$-modules, sont isomorphes. Pour cela nous utilisons les constructions de (\cite{Lu}) dans le cas des big\`{e}bres de Lie bas\'{e}es sur la th\'{e}orie des alg\`{e}bres de Clifford (\cite{K-S}). Nous avons le r\'{e}sultat suivant :
\begin{thm} Soit $(\mathcal{G}, \mu, \gamma, \phi)$ une quasi-big\`{e}bre de Lie. L'application lin\'{e}aire
$$\Re : \mathcal{D}\rightarrow End(\Lambda\mathcal{G}) : x + \xi \rightarrow \Re_{x} + \Re_{\xi}$$
d\'{e}finie par
$$\Re_{x}(Y) = d_{\gamma}(x)\wedge Y + ad_{x}^{\mu}Y - \frac{1}{2}<\xi^{\mu}, x>Y$$
$$\Re_{\xi}(Y) = - i_{d_{\mu}(\xi)}Y + ad_{\xi}^{\gamma *}Y - (i_{\xi}\phi)\wedge Y + \frac{1}{2}<x^{\gamma}, \xi>Y$$
pour $x \in  \mathcal{G}$,  $\xi\in \mathcal{G}^{*}$ et $Y \in \Lambda\mathcal{G}$, est une repr\'{e}sentation de  $\mathcal{D}$ sur $\Lambda\mathcal{G}$.
\end{thm}
{\bf D\'{e}monstration :} Pour montrer que $\Re$ est une repr\'{e}sentation de  $\mathcal{D}$ sur $\Lambda\mathcal{G}$, il suffit d'\'{e}tablir les relations suivantes :
$$\Re_{[x, y]_{\mathcal{D}}} = \Re_{\mu(x, y)} = [\Re_{x}, \Re_{y}], \forall x, y \in  \mathcal{G};$$
$$\Re_{[x, \xi]_{\mathcal{D}}} = - \Re_{ad_{\xi}^{\gamma *}x} + \Re_{ad_{x}^{\mu *}\xi} = [\Re_{x}, \Re_{\xi}], \forall x \in  \mathcal{G}, \forall \xi\in \mathcal{G}^{*};$$
$$\Re_{[\xi, \eta]_{\mathcal{D}}} = \Re_{\phi(\xi, \eta)} + \Re_{\gamma(\xi, \eta)} = [\Re_{\xi}, \Re_{\eta}], \forall \xi, \eta \in \mathcal{G}^{*}.$$
Prouvons tout d'abord que $\Re_{\mu(x, y)} = [\Re_{x}, \Re_{y}], \forall x, y \in  \mathcal{G};$ en effet, par d\'{e}finition, puis en utilisant la condition de 1-cocycle et la relation (3.21) du lemme 3.3, on a pour tout $Y \in \Lambda\mathcal{G}$ :
$$\begin{array}{ccc}
\Re_{\mu(x, y)}(Y) & = &  d_{\gamma}(\mu(x, y))\wedge Y + ad_{\mu(x, y)}^{\mu}Y - \frac{1}{2}<\xi^{\mu}, \mu(x, y)>Y\\
&   & \\
& = & [d_{\gamma}(x), y]^{\mu}\wedge Y +  [x, d_{\gamma}(y)]^{\mu}\wedge Y  + [ad_{x}^{\mu}, ad_{y}^{\mu}](Y).
\end{array}
$$
Par ailleurs
$$\begin{array}{ccc}
[\Re_{x}, \Re_{y}](Y) & = & \Re_{x}(\Re_{y}(Y)) - \Re_{y}(\Re_{x}(Y))\\
 &   & \\
 &  = & \Re_{x}(d_{\gamma}(y)\wedge Y + ad_{y}^{\mu}Y - \frac{1}{2}<\xi^{\mu}, y>Y)\\
 &   &  \\
 &   &  - \Re_{y}(d_{\gamma}(x)\wedge Y + ad_{x}^{\mu}Y - \frac{1}{2}<\xi^{\mu}, x>Y)\\
  &   & \\
 &  = & [x, d_{\gamma}(y)\wedge Y]^{\mu} + d_{\gamma}(x)\wedge ad_{y}^{\mu}Y  -[y, d_{\gamma}(x)\wedge Y]^{\mu}\\
  &   & \\
 &   & - d_{\gamma}(y)\wedge ad_{x}^{\mu}Y + [ad_{x}^{\mu}, ad_{y}^{\mu}](Y)\\
  &   & \\
 &  = & [d_{\gamma}(x), y]^{\mu}\wedge Y + [x, d_{\gamma}(y)]^{\mu}\wedge Y + [ad_{x}^{\mu}, ad_{y}^{\mu}](Y).
\end{array}
$$
D'o\`{u} $\Re_{\mu(x, y)} = [\Re_{x}, \Re_{y}], \forall x, y \in  \mathcal{G}.$\\
Prouvons a pr\'{e}sent que $- \Re_{ad_{\xi}^{\gamma *}x} + \Re_{ad_{x}^{\mu *}\xi} = [\Re_{x}, \Re_{\xi}], \forall x \in  \mathcal{G}, \forall \xi\in \mathcal{G}^{*};$ en effet, par d\'{e}finition, puis en utilisant la relation (3.14) du lemme 3.1 et la proposition 3.2, on a pour tout $Y \in \Lambda\mathcal{G}$ :
$$\begin{array}{ccc}
- \Re_{ad_{\xi}^{\gamma *}x}(Y) + \Re_{ad_{x}^{\mu *}\xi}(Y) & = & -  (d_{\gamma}(ad_{\xi}^{\gamma *}x))\wedge Y - ad_{ad_{\xi}^{\gamma *}x}^{\mu}Y \\
 &   &\\
 &   & + \frac{1}{2}<\xi^{\mu}, ad_{\xi}^{\gamma *}x>Y - i_{d_{\mu}(ad_{x}^{\mu *}\xi)}Y \\
 &   & \\
 &   & + ad_{ad_{x}^{\mu *}\xi}^{\gamma *}Y - (i_{ad_{x}^{\mu *}\xi}\phi)\wedge Y \\
 &   &\\
 &   & + \frac{1}{2}<x^{\gamma}, ad_{x}^{\mu *}\xi>Y.\\
 &   &\\
 & = & - (ad_{\xi}^{\gamma *}d_{\gamma}(x))\wedge Y - (ad_{x}^{\mu}i_{\xi}\phi)\wedge Y \\
 &    &\\
 &    &  - [ad_{x}^{\mu}, i_{d_{\mu}\xi}](Y) + [ad_{x}^{\mu}, ad_{\xi}^{\gamma *}](Y)\\
 &    &\\
 &    & + \sum_{i}(-1)^{i+1}ad_{ad_{y_{i}}^{\mu *}\xi}^{\gamma *}x\wedge\hat{Y_{i}} + (i_{d_{\mu}(\xi)}d_{\gamma}(x)).
\end{array}$$
En utilisant le lemme 3.2, on obtient
$$\begin{array}{ccc}
 - \Re_{ad_{\xi}^{\gamma *}x}(Y) + \Re_{ad_{x}^{\mu *}\xi}(Y) &  = & - (ad_{\xi}^{\gamma *}d_{\gamma}(x))\wedge Y - (ad_{x}^{\mu}i_{\xi}\phi)\wedge Y \\
  &    &\\
  &    &  - [ad_{x}^{\mu}, i_{d_{\mu}\xi}](Y) + [ad_{x}^{\mu}, ad_{\xi}^{\gamma *}](Y)\\
  &    &\\
  &    & + i_{d_{\mu}(\xi)}(d_{\gamma}(x)\wedge Y) - (d_{\gamma}(x))\wedge i_{d_{\mu}(\xi)}Y.
  \end{array}
  $$
 D'autre part
 $$\begin{array}{ccc}
 [\Re_{x}, \Re_{\xi}](Y) & = & \Re_{x}(\Re_{\xi}(Y)) - \Re_{\xi}(\Re_{x}(Y))\\
 &   &\\
  &  = & \Re_{x}(- i_{d_{\mu}(\xi)}Y + ad_{\xi}^{\gamma *}Y - (i_{\xi}\phi)\wedge Y + \frac{1}{2}<x^{\gamma}, \xi>Y)\\
  &   &\\
  &   & \Re_{\xi}(d_{\gamma}(x)\wedge Y + ad_{x}^{\mu}Y - \frac{1}{2}<\xi^{\mu}, x>Y)\\
  &    &\\
  & = & - (ad_{\xi}^{\gamma *}d_{\gamma}(x))\wedge Y - (ad_{x}^{\mu}i_{\xi}\phi)\wedge Y \\
  &    &\\
  &    &  - [ad_{x}^{\mu}, i_{d_{\mu}(\xi)}](Y) + [ad_{x}^{\mu}, ad_{\xi}^{\gamma *}](Y)\\
  &    &\\
  &   & + i_{d_{\mu}(\xi)}(d_{\gamma}(x)\wedge Y) - d_{\gamma}(x)\wedge i_{d_{\mu}(\xi)}Y.
 \end{array}
 $$
 Par comparaison, on trouve que $- \Re_{ad_{\xi}^{\gamma *}x} + \Re_{ad_{x}^{\mu *}\xi} = [\Re_{x}, \Re_{\xi}], \forall x \in  \mathcal{G}, \forall \xi\in \mathcal{G}^{*}.$\\
 Montrons enfin que $\Re_{\phi(\xi, \eta)} + \Re_{\gamma(\xi, \eta)} = [\Re_{\xi}, \Re_{\eta}], \forall \xi, \eta \in \mathcal{G}^{*};$ en effet, par d\'{e}finition et en utilisant la condition de 1-cocycle, la relation 3.22 du lemme 3.3 et la proposition 3.2, on a pour tout $Y \in \Lambda\mathcal{G}$ :
$$\begin{array}{ccc}
 \Re_{\phi(\xi, \eta)}(Y) + \Re_{\gamma(\xi, \eta)}(Y) & = & (d_{\gamma}(\phi(\xi, \eta)))\wedge Y +  ad_{\phi(\xi, \eta)}^{\mu}Y\\
 &  &\\
 &   & - \frac{1}{2}<\xi^{\mu}, \phi(\xi, \eta)>Y - i_{d_{\mu}(\gamma(\xi, \eta))}(Y) + ad_{\gamma(\xi, \eta)}^{\gamma *}Y\\
 &   &\\
 &  & -(i_{\gamma(\xi, \eta)}\phi)\wedge Y + \frac{1}{2}<x^{\gamma}, \gamma(\xi, \eta)>Y\\
 &  & \\
 & = &(d_{\gamma}(\phi(\xi, \eta)))\wedge Y - i_{[d_{\mu}(\xi), \eta]^{\gamma}}(Y) - i_{[\xi, d_{\mu}(\eta)]^{\gamma}}(Y) \\
 &   &  \\
 &   &  + [ad_{\xi}^{\gamma *}, ad_{\eta}^{\gamma *}](Y)-(i_{\gamma(\xi, \eta)}\phi)\wedge Y\\
 &   & \\
 &   & + \sum_{i=1}^{m}(-1)^{i}(\phi(ad_{y_{i}}^{\mu *}\xi, \eta) + \phi(\xi, ad_{y_{i}}^{\mu *}\eta))\wedge\hat{Y_{i}}.
 \end{array}
 $$
Par ailleurs, en utilisant (3.15) et la proposition 2.2 adapt\'{e}e \`{a} $[ , ]^{\gamma}$, on a
$$\begin{array}{ccc}
[\Re_{\xi}, \Re_{\eta}](Y) & = & \Re_{\xi}(- i_{d_{\mu}(\eta)}Y + ad_{\eta}^{\gamma *}Y - (i_{\eta}\phi)\wedge Y + \frac{1}{2}<x^{\gamma}, \eta>Y)\\
&  &\\
 &   & - \Re_{\eta}(- i_{d_{\mu}(\xi)}Y + ad_{\xi}^{\gamma *}Y - (i_{\xi}\phi)\wedge Y + \frac{1}{2}<x^{\gamma}, \xi>Y)\\
 &   & \\
 & =  & - (ad_{\xi}^{\gamma *}i_{d_{\mu}(\eta)} - i_{d_{\mu}(\eta)}ad_{\xi}^{\gamma *})(Y)\\
 &    &\\
 &   & + (ad_{\eta}^{\gamma *}i_{d_{\mu}(\xi)} - i_{d_{\mu}(\xi)}ad_{\eta}^{\gamma *})(Y)\\
 &   &\\
 &   & + [ad_{\xi}^{\gamma *}, ad_{\eta}^{\gamma *}](Y) + (i_{d_{\mu}(\xi)}i_{\eta}\phi)\wedge Y\\
 &    &\\
 &   & - (i_{d_{\mu}(\eta)}i_{\xi}\phi)\wedge Y - (ad_{\xi}^{\gamma *}i_{\eta}\phi)\wedge Y + (ad_{\eta}^{\gamma *}i_{\xi}\phi)\wedge Y \\
 &   &\\
 &  = & (d_{\gamma}(\phi(\xi, \eta)))\wedge Y - i_{[d_{\mu}(\xi), \eta]^{\gamma}}(Y) - i_{[\xi, d_{\mu}(\eta)]^{\gamma}}(Y)\\
 &   &\\
 &   & + [ad_{\xi}^{\gamma *}, ad_{\eta}^{\gamma *}](Y) -(i_{\gamma(\xi, \eta)}\phi)\wedge Y \\
 &   &\\
 &   & + i_{d_{\mu}(\xi)}((i_{\eta}\phi)\wedge Y) - (i_{\eta}\phi)\wedge i_{d_{\mu}(\xi)}Y\\
 &   & \\
 &   & - i_{d_{\mu}(\eta)}((i_{\xi}\phi)\wedge Y) + (i_{\xi}\phi)\wedge i_{d_{\mu}(\eta)}Y.
\end{array}
$$
En comparant les deux expressions, le lemme 3.2 nous permet de conclure que $\Re_{\phi(\xi, \eta)} + \Re_{\gamma(\xi, \eta)} = [\Re_{\xi}, \Re_{\eta}], \forall \xi, \eta \in \mathcal{G}^{*}.$ \\
Ce qui ach\`{e}ve la d\'{e}monstration du th\'{e}or\`{e}me. $\bigtriangleup$\\
{\bf Remarque 4.1 :} Si $\phi = 0$, on retrouve le r\'{e}sultat de Lu (\cite{Lu}) dans le cas des big\`{e}bres de Lie.\\
{\bf Remarque 4.2 :} La repr\'{e}sentation ci-dessus d\'{e}crite ne pr\'{e}serve pas la graduation dans $\Lambda\mathcal{G}$.
\begin{cor} L'application suivante
$$\Gamma : x + \xi \in \mathcal{D}\rightarrow \Gamma_{(x + \xi)} :End(\Lambda\mathcal{G})\rightarrow End(\Lambda\mathcal{G})$$
d\'{e}finie par
$$\Gamma_{(x + \xi)}(T) = \Re_{(x + \xi)}T - T\Re_{(x + \xi)}$$
pour $x \in  \mathcal{G}$,  $\xi\in \mathcal{G}^{*}$ et $T \in End(\Lambda\mathcal{G})$, est une repr\'{e}sentation de  $\mathcal{D}$ sur $End(\Lambda\mathcal{G})$.
\end{cor}

Pour d\'{e}finir l'isomorphisme de $\mathcal{D}$-modules entre $\Lambda\mathcal{D}$ et $End(\Lambda\mathcal{G})$, comme dans le cas des big\`{e}bres de Lie (\cite{Lu}), on introduit l'\'{e}l\'{e}ment
$$ exp_{\wedge}{\bf r} = {\bf r} + \frac{1}{2!}{\bf r}\wedge{\bf r} + \frac{1}{3!}{\bf r}\wedge{\bf r} +... \in \Lambda\mathcal{D}$$
o\`{u}
$$ {\bf r} = \frac{1}{2}\sum_{i}e_{i}\wedge \xi^{i}$$
est l'\'{e}l\'{e}ment d\'{e}finissant la structure canonique de quasi-big\`{e}bre de Lie sur $\mathcal{D}$. Pour $A = \xi_{i}\wedge...\wedge\xi_{k}\in \mathcal{G}^{*}$, posons $\hat{A} = \xi_{k}\wedge...\wedge\xi_{1}$. On a le r\'{e}sultat suivant :
\begin{thm} Pour $U\in \Lambda\mathcal{D}$, posons
$$i_{exp_{\wedge}{\bf r}}U = \sum_{j}X_{j}\otimes A_{j} \in \Lambda\mathcal{D} \cong \Lambda\mathcal{G}\otimes\Lambda\mathcal{G}^{*}$$
o\`{u} $X_{j}\in\Lambda\mathcal{G}$ et $A_{j}\in\Lambda\mathcal{G}^{*}$. D\'{e}finissons
$$ Q(U)\in End(\Lambda\mathcal{G}) : Q(U)(Y) = \sum_{j}X_{j}\wedge i_{\hat{A_{j}}}Y.$$
Alors l'application
$$Q : \Lambda\mathcal{D}\rightarrow End(\Lambda\mathcal{G})$$
est un isomorphisme de $\mathcal{D}$-modules, o\`{u} $\mathcal{D}$ agit sur $\Lambda\mathcal{D}$ par l'action adjointe et sur $End(\Lambda\mathcal{G})$ par l'application $\Gamma$.
\end{thm}
{\bf D\'{e}monstration :} Il s'agit de montrer que $\forall  x + \xi \in  \mathcal{D}$, on a
$$\Gamma_{(x + \xi)}\circ Q = Q\circ ad_{(x + \xi)}.$$
Plus pr\'{e}cisement, $\forall x + \xi \in  \mathcal{D}$, $\forall U\in \Lambda\mathcal{D}$ et $\forall Y\in \Lambda\mathcal{G}$, on doit avoir
$$ (\Gamma_{(x + \xi)}(Q(U)))(Y) = (Q(ad_{(x + \xi)}(U)))(Y).$$
En effet, $\forall x \in  \mathcal{G}$, $\forall U = X \in \Lambda\mathcal{G}$ et $\forall Y\in \Lambda\mathcal{G}$, on a par d\'{e}finition de $Q$, $Q(X))(Y) = X\wedge Y$; d'o\`{u}
$$\begin{array}{ccc}
\Gamma_{x}(Q(X)))(Y) & = & \Re_{x}(Q(X)(Y)) - Q(X)(\Re_{x}(Y))\\
&     &\\
&  =  & \Re_{x}(X\wedge Y) - Q(X)(d_{\gamma}(x)\wedge Y + ad_{x}^{\mu}Y - \frac{1}{2}<\xi^{\mu}, x>Y)\\
&     &\\
&  =  & d_{\gamma}(x)\wedge X\wedge Y + ad_{x}^{\mu}(X\wedge Y) - \frac{1}{2}<\xi^{\mu}, x>(X\wedge Y)\\
&     &\\
&     & -X\wedge d_{\gamma}(x)\wedge Y - X\wedge (ad_{x}^{\mu}Y) + \frac{1}{2}<\xi^{\mu}, x>(X\wedge Y)\\
&     &\\
&  =  & (ad_{x}X)\wedge Y = Q(ad_{x}(X))(Y).
\end{array}$$
Ce qui prouve la relation pour $ x \in  \mathcal{G}$ et $U = X \in \Lambda\mathcal{G}.$\\
Prouvons la relation pour $ \xi \in  \mathcal{G}^{*}$ et $U = X \in \Lambda\mathcal{G}.$ En effet, $\forall Y\in \Lambda\mathcal{G}$\\
$$\begin{array}{ccc}
\Gamma_{\xi}(Q(X))(Y) & = & \Re_{\xi}(Q(X)(Y)) - Q(X)(\Re_{\xi}(Y))\\
&     &\\
&  =  & \Re_{\xi}(X\wedge Y) - Q(X)(- i_{d_{\mu}(\xi)}Y + ad_{\xi}^{\gamma ^*}Y  \\
&     &\\
&     & - (i_{\xi}\phi)\wedge Y + \frac{1}{2}<x^{\gamma}, \xi>Y)\\
&     &\\
&  =  & - i_{d_{\mu}(\xi)}(X\wedge Y) + ad_{\xi}^{\gamma ^*}(X\wedge Y)\\
&     &\\
&     & - (i_{\xi}\phi)\wedge (X\wedge Y)+ \frac{1}{2}<x^{\gamma}, \xi>(X\wedge Y)\\
&     &\\
&     & + X\wedge (i_{d_{\mu}(\xi)}Y) -  X\wedge (ad_{\xi}^{\gamma ^*}Y) \\
&     &\\
&     &  + X\wedge(i_{\xi}\phi)\wedge Y - \frac{1}{2}<x^{\gamma}, \xi>(X\wedge Y)\\
&     &\\
&  =  & - i_{d_{\mu}(\xi)}(X\wedge Y) + X\wedge (i_{d_{\mu}(\xi)}Y) + (ad_{\xi}^{\gamma ^*}X)\wedge Y .
\end{array}$$
Pour calculer $Q(ad_{\xi}(X))(Y)$, $\mathcal{G}$ \'{e}tant suppos\'{e}e de dimension finie, on peut consid\'{e}rer $X$ et $Y$ comme des \'{e}l\'{e}ments d\'{e}composables de $\Lambda\mathcal{G}$, i.e $X = x_{1}\wedge ...\wedge x_{m}$ et $Y = y_{1}\wedge ...\wedge y_{n}$; par d\'{e}finition de $[ , ]_{\mathcal{D}}$
$$\begin{array}{ccc}
ad_{\xi}X & = & ad_{\xi}^{\gamma *}X - \sum_{k=1}^{m}x_{1}\wedge ...\wedge (ad_{x_{k}}^{\mu ^{*}}\xi)\wedge...\wedge x_{m}\\
&     & \\
&   =  & ad_{\xi}^{\gamma *}X - \sum_{k=1}^{m}(-1)^{n-k}\hat{X}_{k}\wedge(ad_{x_{k}}^{\mu ^{*}}\xi),
\end{array}
$$
o\`{u} $\hat{X}_{k} = x_{1}\wedge ...\wedge \hat{x_{k}}\wedge...\wedge x_{m}$. Ce qui nous permet d'avoir
$$\begin{array}{ccc}
i_{\bf r}(ad_{\xi}X) &  = & - \frac{1}{2}(-1)^{m-1}\sum_{k=1}^{m}\sum_{i}(-1)^{m-k}i_{\xi^{i}}(\hat{X}_{k})\wedge i_{e_{i}}(ad_{x_{k}}^{\mu ^{*}}\xi)\\
&   &\\
&  = & - \frac{1}{2}\sum_{k=1}^{m}\sum_{i}\sum_{l=1}^{m}(-1)^{k+l}(i_{\xi^{i}}x_{l})\hat{X}_{kl}\wedge (i_{e_{i}}(ad_{x_{k}}^{\mu ^{*}}\xi))\\
&    &\\
&  =  & - \frac{1}{2}\sum_{k=1}^{m}\sum_{l=1}^{m}(-1)^{k+l}(i_{x_{l}}(ad_{x_{k}}^{\mu ^{*}}\xi))\hat{X}_{kl}\\
&   & \\
&  =  & - \sum_{k=1}^{m}\sum_{l=1}^{m}(-1)^{k+l}(i_{d_{\mu}\xi}(x_{k}\wedge x_{l}))\hat{X}_{kl}\\
&   &\\
&  =  & -i_{d_{\mu}\xi}X
\end{array}
$$
o\`{u} $ \hat{X}_{kl} = x_{1}\wedge ...\wedge \hat{x_{k}}\wedge...\wedge\hat{x_{l}}\wedge...\wedge x_{m}$. Ainsi
$$\begin{array}{ccc}
i_{exp_{\wedge}{\bf r}}(ad_{\xi}X) & = & ad_{\xi}X + i_{\bf r}(ad_{\xi}X)\\
&   &\\
&  =  & ad_{\xi}^{\gamma *}X - \sum_{k=1}^{m}(-1)^{m-k}\hat{X_{k}}\wedge(ad_{x_{k}}^{\mu ^{*}}\xi) - i_{d_{\mu}\xi}X.
\end{array}$$
Par cons\'{e}quent
$$\begin{array}{ccc}
Q(ad_{\xi}(X))(Y) & = & (ad_{\xi}^{\gamma *}X)\wedge Y - (i_{d_{\mu}\xi}X)\wedge Y \\
 &  &\\
 &   & + \sum_{k=1}^{m}(-1)^{m-k}\hat{X_{k}}\wedge(i_{ad_{x_{k}}^{\mu ^{*}}\xi}Y)\\
 &   &\\
 &  = & (ad_{\xi}^{\gamma *}X)\wedge Y - (i_{d_{\mu}\xi}X)\wedge Y \\
 &    & \\
 &    & + \sum_{k=1}^{m}\sum_{l=1}^{n}(-1)^{m-k-l}(i_{d_{\mu}\xi}(x_{k}\wedge y_{l}))\hat{X_{k}}\wedge \hat{Y_{l}}.
\end{array}$$
Mais
$$\begin{array}{ccc}
i_{d_{\mu}\xi}(X\wedge Y) & = & (i_{d_{\mu}\xi}X)\wedge Y  + X\wedge(i_{d_{\mu}\xi}X)\\
 &   &\\
 &   & - \sum_{k=1}^{m}\sum_{l=1}^{n}(-1)^{m-k-l}(i_{d_{\mu}\xi}(x_{k}\wedge y_{l}))\hat{X_{k}}\wedge \hat{Y_{l}}.
\end{array}
$$
En d\'{e}finitive, on trouve
$$Q(ad_{\xi}(X))(Y)  = - i_{d_{\mu}(\xi)}(X\wedge Y) + X\wedge (i_{d_{\mu}(\xi)}Y) + (ad_{\xi}^{\gamma ^*}X)\wedge Y .$$
La relation est donc prouv\'{e}e pour tout $x + \xi \in  \mathcal{D}$, $\forall U = X\in \Lambda\mathcal{G}$ et $\forall Y\in \Lambda\mathcal{G}.$\\
Pour tout $x \in  \mathcal{G}$, $\forall U = A \in \Lambda\mathcal{G}^{*}$ et $\forall Y\in \Lambda\mathcal{G}$, on a par d\'{e}finition de $Q$
$$Q(A)(Y) = i_{\hat{A}}Y$$
et
$$\begin{array}{ccc}
\Gamma_{x}(Q(A))(Y) & = & \Re_{x}(Q(A)(Y)) - Q(A)(\Re_{x}(Y))\\
&  &\\
& = & \Re_{x}(i_{\hat{A}}Y) - i_{\hat{A}}(\Re_{x}(Y))\\
&   &\\
& = & d_{\gamma}(x)\wedge i_{\hat{A}}Y + ad_{x}^{\mu}i_{\hat{A}}Y - \frac{1}{2}<\xi^{\mu}, x>i_{\hat{A}}Y\\
&   &\\
&   & - i_{\hat{A}}(d_{\gamma}(x)\wedge Y) - i_{\hat{A}}ad_{x}^{\mu}Y + \frac{1}{2}<\xi^{\mu}, x>i_{\hat{A}}Y\\
&   &\\
&  = & - i_{\hat{A}}(d_{\gamma}(x)\wedge Y) + d_{\gamma}(x)\wedge i_{\hat{A}}Y + [ad_{x}^{\mu}, i_{\hat{A}}](Y).
\end{array}$$
Calculons $Q(ad_{x}(A))(Y)$; pour cela supposons que $A = \xi_{1}\wedge ...\wedge \xi_{m}$ avec $m \leq |Y|$. Par d\'{e}finition de $[ , ]_{\mathcal{D}}$
$$\begin{array}{ccc}
ad_{x}A & = & \sum_{k =1}^{m}(-1)^{k}(ad_{\xi_{k}}^{\gamma ^{*}}x)\wedge \hat{A}_{k} + ad_{x}^{\mu *}A\\
 &    & \\
 &  =  & \sum_{k =1}^{m}(-1)^{k}(i_{\xi_{k}}d_{\gamma}(x))\wedge \hat{A}_{k} + ad_{x}^{\mu *}A
\end{array}$$
d'o\`{u}
$$\begin{array}{ccc}
i_{exp_{\wedge}{\bf r}}(ad_{x}A) & = & ad_{x}A + i_{\bf r}(ad_{x}A)\\
&   &\\
& = & \sum_{k =1}^{m}(-1)^{k}(i_{\xi_{k}}d_{\gamma}(x))\wedge \hat{A}_{k} + ad_{x}^{\mu *}A + i_{d_{\gamma}(x)}A.
\end{array}$$
Ainsi, par d\'{e}finition de $Q$
$$Q(ad_{x}(A))(Y) = i_{ad_{x}^{\mu *}\hat{A}}Y + \sum_{k =1}^{m}(-1)^{k}(i_{\xi_{k}}d_{\gamma}(x))\wedge i_{\hat{A}_{k}}Y + i_{i_{d_{\gamma}(x)}\hat{A}}Y.$$
Mais
$$i_{ad_{x}^{\mu *}\hat{A}}Y = [ad_{x}^{\mu}, i_{\hat{A}}](Y);$$
$$i_{i_{d_{\gamma}(x)}\hat{A}}Y = \sum_{k,l=1}^{m}(-1)^{k+l}(i_{\xi_{l}}i_{\xi_{k}}d_{\gamma}(x))i_{\hat{A}_{kl}}Y$$
et
$$\begin{array}{ccc}
i_{\hat{A}}(d_{\gamma}(x)\wedge Y) & = & - \sum_{k=1}^{m}\sum_{l=1}^{m}(-1)^{k+l}(i_{\xi_{l}}i_{\xi_{k}}d_{\gamma}(x))i_{\hat{A}_{kl}}Y\\
&   &\\
&   & + \sum_{k=1}^{m}(-1)^{k-1}(i_{\xi_{k}}d_{\gamma}(x))i_{\hat{A}_{k}}Y + d_{\gamma}(x)\wedge i_{\hat{A}}Y.
\end{array}$$
En d\'{e}finitive
$$Q(ad_{x}(A))(Y) = - i_{\hat{A}}(d_{\gamma}(x)\wedge Y) + d_{\gamma}(x)\wedge i_{\hat{A}}Y + [ad_{x}^{\mu}, i_{\hat{A}}](Y).$$
D'o\`{u}
$$\Gamma_{x}(Q(A))(Y) = Q(ad_{x}(A))(Y), \hspace{0.5cm} \forall x \in \mathcal{G}, \forall U = A \in \Lambda\mathcal{G}^{*},  \forall Y\in \Lambda\mathcal{G}.$$
Pour tout $\xi \in  \mathcal{G}^{*}$, $\forall U = A \in \Lambda\mathcal{G}^{*}$ et $\forall Y\in \Lambda\mathcal{G}$, on a par d\'{e}finition de $Q$
$$Q(A)(Y) = i_{\hat{A}}Y$$
et
$$\begin{array}{ccc}
\Gamma_{\xi}(Q(A))(Y) & = & \Re_{\xi}(Q(A)(Y)) - Q(A)(\Re_{\xi}(Y))\\
&  &\\
& = & \Re_{\xi}(i_{\hat{A}}Y) - i_{\hat{A}}(\Re_{\xi}(Y))\\
&   &\\
&  = & [ad_{\xi}^{\gamma *}, i_{\hat{A}}](Y) -(i_{\xi}\phi)\wedge i_{\hat{A}}Y + i_{\hat{A}}((i_{\xi}\phi)\wedge Y)\\
&    &\\
&  = & [ad_{\xi}^{\gamma *}, i_{\hat{A}}](Y) + [i_{\hat{A}}, \varepsilon_{i_{\xi}\phi}](Y).
\end{array}$$
En supposant que $A$ est un \'{e}l\'{e}ment d\'{e}composable de $\Lambda\mathcal{G}^{*}$, i.e $A = \xi_{1}\wedge ...\wedge \xi_{m}$, on a par d\'{e}finition de $[ , ]_{\mathcal{D}}$
$$ad_{\xi}A = ad_{\xi}^{\gamma}A + \sum_{k=1}^{m}(-1)^{k-1}\phi(\xi, \xi_{k})\wedge\hat{A}_{k}.$$
Par un calcul, on trouve
$$Q(ad_{\xi}(A))(Y) = [ad_{\xi}^{\gamma *}, i_{\hat{A}}](Y) + [i_{\hat{A}}, \varepsilon_{i_{\xi}\phi}](Y).$$
Par cons\'{e}quent
$$\Gamma_{(x + \xi)}(Q(A))(Y) = Q(ad_{(x + \xi)}(A))(Y), \hspace{0.5cm} \forall x + \xi \in \mathcal{D}, \forall U = A \in \Lambda\mathcal{G}^{*},  \forall Y\in \Lambda\mathcal{G}.$$
Pour une d\'{e}monstration compl\`{e}te du th\'{e}or\`{e}me dans le cas g\'{e}n\'{e}ral, nous allons nous servir d'un r\'{e}sultat de (\cite{Lu}) o\`{u} la d\'{e}monstration est faite pour le cas des big\`{e}bres de Lie. Pour cela, consid\'{e}rons l'op\'{e}rateur
$$D : \Lambda\mathcal{G}\otimes\Lambda\mathcal{G}^{*} \rightarrow End(\Lambda\mathcal{G}) : D_{X\wedge A}(Y) = X\wedge i_{\hat{A}}Y.$$
Alors d'apr\`{e}s (\cite{Lu}), l'application
$$D\circ i_{exp_{\wedge}{\bf r}} : \Lambda\mathcal{G}\otimes\Lambda\mathcal{G}^{*} \rightarrow End(\Lambda\mathcal{G})$$
est un isomorphisme de $\mathcal{D}$-modules d'alg\`{e}bre de Lie. Pour conclure la d\'{e}monstration du th\'{e}or\`{e}me, il suffit de remarquer que $Q = D\circ i_{exp_{\wedge}{\bf r}}$. $\bigtriangleup$

\vspace{1cm}
{\bf Remerciements :} L'auteur remercie vivement le Centre Abdus Salam ICTP pour l'hospitalit\'{e} et le soutien mat\'{e}riel qui ont permis la r\'{e}alisation de ce travail.

This work was done within the Associateship Scheme of the Abdus Salam International Centre for Theoretical Physics, Trieste, Italy. Financial support from ICTP is acknowledged.

\end{document}